\input amstex
\magnification=1200
\documentstyle{amsppt}
\NoBlackBoxes
\rightheadtext{Harish-Chandra modules with arbitrary minimal $\k$-type}
\widestnumber\key{[Most]2]}
\baselineskip 20pt

\hsize 6.25truein
\vsize = 8.75truein


\catcode`\@=11
\def\logo@{}
\catcode`\@=13

\def\cplus{\hbox{$\subset${\raise0.3ex\hbox{\kern -0.55em
${\scriptscriptstyle +}$}}\ }}

\def\clplus{\hbox{$\subset${\raise0.3ex\hbox{\kern -0.55em ${\scriptscriptstyle +}$}}\ }}
\def\crplus{\hbox{$\supset${\raise1.05pt\hbox{\kern -0.55em ${\scriptscriptstyle +}$}}\ }}

\def\N{\bold N}

\def\p{\bold p}
\def\C{\bold C}
\def\Z{\bold Z}
\def\R{\bold R}

\def\supp{\text{\rm supp}}
\def\Re{\text{\rm Re}}

\def\s{\bold S}

\def\p{\bold p}
\def\C{\bold C}
\def\Z{\bold Z}
\def\R{\bold R}

\def\ch{\text{\rm ch}}

\def\Hom{\text{\rm Hom}}
\def\m{\frak m}
\def\n{\frak n}
\def\p{\frak p}
\def\g{\frak g}
\def\h{\frak h}
\def\b{\frak b}
\def\k{\frak k}
\def\t{\frak t}

\def\s{\frak s}

\vskip .15in
\topmatter
\title A CONSTRUCTION OF GENERALIZED HARISH-CHANDRA MODULES WITH ARBITRARY MINIMAL $\frak k$-type \endtitle
\author Ivan Penkov and Gregg Zuckerman \endauthor
\affil  \endaffil
\abstract  Let $\frak g$ be a semisimple complex Lie algebra and $\k\subset\g$ be any algebraic subalgebra reductive in $\frak g$.  For any simple finite dimensional $\frak k$-module $V$, we construct simple $(\frak g,\frak k)$-modules $M$ with finite dimensional $\frak k$-isotypic components such that $V$ is a $\frak k$-submodule of $M$ and the Vogan norm of any simple $\k$-submodule $V'\subset M, V'\not\simeq V$, is greater than the Vogan norm of $V$.  The $(\frak g,\frak k)$-modules $M$ are subquotients of the fundamental series of $(\frak g,\frak k)$-modules introduced in \cite{PZ2}.
\endabstract
\endtopmatter

\head Introduction \endhead

The structure theory of infinite dimensional modules over finite dimensional semisimple Lie algebras has its roots in the description of all finite dimensional representations. Celebrated landmarks of the theory are the classification of simple Harish-Chandra modules and the computation of the characters of simple highest weight modules (the Kazhdan-Lusztig conjecture).  A deep open problem in the structure theory of modules over  a complex semisimple Lie algebra $\frak g$ is the construction and eventual classification of all simple generalized Harish-Chandra modules, see \cite{PZ1}.
By definition, a simple $\g$-module $M$ is  a {\it generalized Harish-Chandra module} if $M$ has finite dimensional isotypic components as module over some reductive in $\frak g$ subalgebra of $\g$.  Equivalently, a simple generalized Harish-Chandra module is a simple $\g$-module $M$ for which the multiplicities of $M$ as a $\frak g [M]$-module are finite. The subalgebra $\g [M]\subset \g$ is defined as the set of all elements of $\g$ which act locally finitely on $M$, see \cite{F} and \cite{PZ1}.  In \cite{PSZ} we have proved that, if the multiplicities of $M$ as a $\g [M]$-module are finite, then $\g [M]$ has a natural reductive part $\g [M]_{\text{red}}$, and that $M$ has finite type also as a $\g [M]_{\text{red}}$-module, i.e. the dimensions of all $\g [M]_{\text{red}}$-isotypic components of $M$ are finite.

Recently two considerable steps in the study of simple generalized Harish-Chandra modules have been made.  In \cite{PSZ} we have described explicitly all possible subalgebras $\g [M]_{\text{red}}\subset\g$ arising from simple generalized Harish-Chandra modules (these are the primal subalgebras of $\g$, see \cite{PSZ}), and in \cite{PZ2} we have classified all simple generalized Harish-Chandra modules $M$ with generic minimal $\k$-type.  Here $\k$ stands for any algebraic reductive in $\g$ subalgebra $\frak k$ with $\k\subset\g [M]$ such that $M$ has finite dimensional $\frak k$-isotypic components.  The latter result raises a natural question: for a fixed reductive in $\g$ algebraic subalgebra $\frak k$, what are the minimal $\frak k$-types arising from simple $(\g,\k )$-modules of finite type?  In the case when the pair $(\g,\k )$ is symmetric, it is known from Vogan's classification of Harish-Chandra modules that there is no obstruction for a simple finite dimensional $\frak k$-module to be the minimal $\frak k$-type of a simple $(\g,\k )$-module.

The purpose of the present note is to give a simple proof of this fact by a direct construction in the case of an arbitrary algebraic reductive in $\g$ subalgebra $\k\subset\g$.  Our construction is based on the fundamental series of $(\g,\k )$-modules, introduced recently in \cite{PZ2}, and extends the construction of a simple $(\g,\k)$-module with an arbitrary minimal $\frak k$-type given in \cite{PZ1} for the case where $\frak k$ is a principal $s\ell (2)$-subalgebra of $\frak g$.

\head 1. Conventions and Preliminaries \endhead

The ground field is $\C$, and if not explicitly stated otherwise, all vector spaces and Lie algebras are defined over $\C$.  By definition, $\N = \{0,1,2,\ldots\}$.  The sign $\otimes$ denotes tensor product over $\C$.  The superscript $*$ indicates dual space, and $\Lambda^\cdot(\quad)$ and $S^\cdot (\quad)$ denote respectively the exterior and symmetric algebra. By $Z(\frak l)$ we denote the center of a Lie algebra $\frak l$, $U(\frak l)$ stands for the enveloping algebra of $\frak l$, and $H^\cdot (\frak l, M)$ stands for the cohomology of a Lie algebra $\frak l$ with coefficients in an $\frak l$-module $M$.
The sign $\cplus$ stands for semidirect sum of Lie algebras (if $\frak l = \frak l'\cplus\frak l''$, then $\frak l'$ is an ideal in $\frak l$ and $\frak l''\simeq\frak l/\frak l'$).

If $\frak l$ is a Lie algebra, $M$ is an $\frak l$-module, and $\omega\in\frak l^*$, we put \newline $M^\omega := \{m\in M\vert\ell\cdot m = \omega (\ell)m \,\,  \forall\ell\in\frak l\}$.  We call $M^\omega$ a {\it weight space} of $M$ and we say that $M$ is an $\frak l$-{\it weight module} if
$$
M = \bigoplus\limits_{\omega\in\frak l^{*}} M^\omega.
$$
By $\supp_{\frak l} M$ we denote the set $\{\omega\in\frak l^*\vert M^\omega\ne 0\}$.

A finite {\it multiset} is a function $f$ from a finite set $D$ into $\N$.  A {\it submultiset} of $f$ is a multiset $f'$ defined on the same domain $D$ such that $f'(d)\leq f(d)$ for any $d\in D$.  For any finite multiset $f$, defined on an additive monoid $D$, we can put $\rho_f := \frac{1}{2}\sum\limits_{d\in D} f(d)d$.  If $M$ is an $\frak l$-weight module as above, and $\dim M < \infty$, $M$ determines the finite multiset  $\ch_{\frak l} M$ which is the function $\omega\mapsto \dim M^\omega$ defined on $\supp_{\frak l} M$.

Let $\frak g$ be a fixed finite dimensional semisimple Lie algebra and $\frak k\subset\frak g$ be a fixed algebraic subalgebra which is reductive in $\frak g$.  Fix a Cartan subalgebra $\frak t$ of $\frak k$ and a Cartan subalgebra $\frak h$ of $\frak g$ such that $\frak t\subset\frak h$.  Note that, since $\frak k$ is reductive in $\frak g$, $\frak g$ is a $\frak t$-weight module. Note also that the $\R$-span of the roots $\Delta$ of $\frak h$ in $\frak g$ fixes a real structure on $\frak h^*$, whose projection onto $\frak t^*$ is a well-defined real structure on $\frak t^*$.  In what follows, we will denote by $\Re \lambda$ the real part of an element $\lambda\in\frak t^*$.  We fix also a Borel subalgebra $\frak b_{\frak k} \subset\frak k$ with $\frak b_{\frak k}\supset\frak t$.  Then $\frak b_{\frak k} = \frak t ~\crplus \frak n_{\frak k}$, where $\frak n_{\frak k}$ is the nilradical of $\frak b_{\frak k}$.   We set $\rho :=\rho_{\ch_{\frak t} \frak n_{\frak k}}$, and by $W_{\frak k}$ we denote the Weyl group of $\frak k$.

Let $ \langle \, ,\, \rangle$ denote the unique $\frak g$-invariant symmetric bilinear form on $\frak g^*$ such that $\langle \alpha,\alpha \rangle =2$ for any long root of a simple component of $\frak g$.  The form $\langle \, , \, \rangle$ enables us to identify $\frak g$ with $\frak g^*$. Then $\frak h$ is identified with $\frak h^*$, and $\frak k$ is identified with $\frak k^*$.  We will sometimes consider $\langle \, ,\, \rangle$ as a form on $\frak g$.  The superscript $\perp$ indicates orthogonal space.  Note that there is a canonical $\k$-module decomposition $\frak g = \frak k\oplus\frak k^\perp$.  We also set $\Vert\kappa\Vert^2 := \langle \kappa,\kappa \rangle$ for any $\kappa\in\frak h^*$.

To any $\lambda\in\frak t^*$ we associate the following parabolic subalgebra $\frak p_\lambda$ of $\frak g$:
$$
\frak p_\lambda = \frak  h \oplus (\bigoplus\limits_{\alpha\in\Delta_{\lambda}} \frak g^\alpha),
$$
where $\Delta_\lambda :=\{\alpha\in\Delta\mid \langle \Re \lambda, \alpha \rangle\geq 0\}$.  By $\frak m_\lambda$ and $\frak n_\lambda$ we denote respectively the reductive part of $\frak p_\lambda$ (containing $\frak h$) and the nilradical of $\frak p_\lambda$.  In particular $\frak p_\lambda = \m_\lambda ~\crplus \n_\lambda$, and if $\lambda$ is $\frak b_{\frak k}$-dominant, then $\frak p_\lambda\cap \frak k = \frak b_{\frak  k}$.    We call $\frak p_\lambda$ a {\it compatible parabolic subalgebra}.  A compatible parabolic subalgebra $\frak p = \frak m ~\crplus \frak n$ (i.e. $\frak p = \frak p_\lambda$ for some $\lambda\in\frak t^*$) is {\it minimal} if it does not properly contain another compatible parabolic subalgebra.  It is an important observation that if $\frak p = \frak m ~\crplus\frak n$ is minimal, then $\frak t\subset Z(\frak m)$.

A {\it $\k$-type} is by definition a simple finite dimensional $\k$-module.  By $V(\mu)$ we will denote a $\frak k$-type with $\frak b_{\frak k}$-highest weight $\mu$ ($\mu$ is then $\frak k$-integral and $\frak b_{\frak k}$-dominant).

For the purposes of this paper, we call a $\g$-module $M$ a $(\g,\k)$-{\it module} if $M$ is isomorphic as a $\k$-module to a direct sum of isotypic components of $\k$-types.  We say that a $(\g,\k)$-module $M$ is {\it of finite type} if $\dim \text{Hom}_\k (V(\mu),M)< \infty$ for every $\k$-type $V(\mu)$. We say also that a $\k$-type $V$ {\it is a $\k$-type of} $M$  if $\dim_\k \text{Hom}(V,M)\neq 0$.
If $M$ is a $(\g,\k )$-module, a $\k$-type $V(\mu)$ of $M$ is {\it minimal} if the Vogan norm, i.e. the function $\mu'\mapsto \Vert\Re \mu' + 2\rho\Vert^2$, defined on the $\b_\k$-highest weights $\mu'$ of all $\k$-types of $M$, has a minimum at $\mu$.  Any simple $(\g,\k)$-module $M$ has a minimal $\k$-type.
\vskip .15in

Recall that the functor of $\k$-locally finite vectors $\Gamma_{\k,\t}$ is a well-defined left exact functor on the category of $(\g,\t)$-modules with values in $(\g,\k)$-modules,
$$
\Gamma_{\k,\t} (M) = \sum\limits_{M'\subset M,\dim M'=1, \dim U(\k)\cdot M' < \infty} M'.
$$
By $R^\cdot \Gamma_{\k,\t}:=\bigoplus\limits_{i\geq 0} R^i \Gamma_{\k,\t}$ we denote as usual the total right derived functor of $\Gamma_{\k,\t}$, see \cite{PZ1} and the references therein.

Let $\p = \m ~\crplus\n$ be a minimal compatible parabolic subalgebra, $E$ be a simple finite dimensional $\p$-module,
$\rho_\n := \rho_{\ch_{\t}\n}$ and $\rho^\perp_{\n} := \rho_{\ch_{\t}(\n\cap\k^{\perp})}$. Set
$$
F^\cdot(\p, E) := R^\cdot\Gamma_{\k,\t} (\Gamma_{\t,0} (\Hom_{U(\p)}
(U(\g), E\otimes\Lambda^{\dim \n}(\n)))).
$$
By definition, $F^\cdot (\p, E)$ is the {\it fundamental series of} $(\g,\k )$-{\it modules}.

\head 2. Main Results \endhead

\proclaim{Theorem 1}  Let $V$ be any $\k$-type.  There exists a simple $(\g,\k)$-module of finite type $M$ such that $V$ is the unique minimal $\k$-type of $M$.
\endproclaim

The proof is based on the following construction.  Let $V = V(\mu )$ be a fixed $\k$-type and let $\p = \m\crplus\frak n$ be any minimal compatible parabolic subalgebra of $\g$ which lies in $\p_{\mu +2\rho}$.  Let, in addition, $E$ be any simple finite dimensional $\p$-module on which $\t$ acts via the weight $\mu - 2\rho^\perp_\n$ ($E$ exists since $\t\subset Z(\m ))$.

\proclaim{Theorem 2}  Let $s = \dim\n_\k$.  The $(\g,\k)$-module $F^s(\p, E)$ is of finite type and is non-zero.  In addition, $V$ is the unique minimal $\k$-type of $F^s(\p,E)$ and $\dim\Hom_\k (V, F^s (\p,E))=\dim E$.
\endproclaim

Theorem 2 implies Theorem 1 as a module $M$ whose existence is claimed by Theorem 1 can be constructed as any simple quotient of a $\g$-submodule of $F^s(\p, E)$ generated by the image of any $\k$-module injection $V\to F^s (\p,E)$.

Theorem 2 is a direct corollary of the following five statements: two more general propositions and three lemmas under the assumptions of Theorem 2.

\proclaim{Proposition 1} Let $\p = \m\crplus\n$ be any minimal parabolic subalgebra, $E$ be any simple finite dimensional $\p$-module, and  $V(\delta)$ be a $\k$-type of $F^{s-i}(\p,E)$ for some $i\in\Z$.

a)  There exists $w\in W_\k$ of length $i$ (in particular, $i\in\N$) and a multiset
$$\align
&n_\cdot := \supp_\t (\n\cap\k^\perp) \to \N,\\
&\qquad \beta\mapsto n_\beta
\endalign
$$
such that
$$
\omega = w(\delta +\rho )-\rho -2\rho^\perp_\n - \sum\limits_{\beta} n_\beta \beta ,
$$
where $\omega$ is the weight via which $\t$ acts on $E$.  Furthermore, $\dim\text{Hom}_\k (V(\delta),F^{s-i}(\p, E))$
is bounded by the integer
$$
\dim E \, (\sum\limits_{\ell (w)=i} \dim (S^\cdot(\frak n\cap \frak k^\perp)^{\xi (w)})),
$$
where $\xi (w) = w(\delta +\rho)-\rho -\omega -2\rho^\perp_\n$, and $S^\cdot (\n\cap\k^\perp )$ is considered as a $\k$-weight module.
\endproclaim

\proclaim{Proposition 2} Under the assumptions of Proposition 1,
$$
\sum\limits_{0\leq i\leq s} (-1)^i \dim \Hom_\k (V(\delta), F^{s-i} (\p, E)) \tag1
$$
$$
= \sum\limits_{0\leq j\leq s} (-1)^j (\sum\limits^\infty_{m=0} \dim \Hom_\t (H^j (\n\cap\k, V(\delta)), S^m(\n\cap\k^\perp)\otimes E\otimes\Lambda^{\dim (\n\cap\k^{\perp})} (\n\cap\k^\perp))),
$$
and the inner sum on the right hand side of (1) is finite.
\endproclaim

Propositions 1 and 2 are a modification of Theorem 6.3.12 and
Corollary 6.3.13 in \cite{V}, and their proofs follow exactly the
same lines (an inspection of Vogan's proofs reveals that the symmetry
assumption on $(\g,\k)$ is not needed).  Therefore, we refer the
reader to \cite{V}.

Proposition 1 implies that, for any minimal compatible parabolic subalgebra $\p$ and for any simple finite dimensional $\p$-module $E$, $F^\cdot (\p, E)$ (and thus $F^s(\p, E)$) is a $(\g,\k)$-module of finite type, and also that $F^i(\p,E) = 0$ for $i > s$.

In the rest of this section we assume that $\p$ and $E$ are as in Theorem 2.

\proclaim{Lemma 1}  If $V = V(\mu )$  is a $\k$-type of $F^{s-i} (\p,E)$, then $i=0$.
\endproclaim

\demo{Proof}  Choose $\lambda\in\h^*$ so that $\p = \p_\lambda$.  In particular, $\langle \Re\lambda,\gamma\rangle > 0$ for $\gamma\in\supp_\t \n$.  By Proposition 1, there exist $w\in W_\k$ of length $i$ and a multiset $n$.: $\supp_\t (\n\cap \k^\perp)\to\N$ such that
$$
\omega = w(\mu + \rho ) - \rho - 2\rho^\perp_\n - \sum\limits_{\beta\in\supp_{\t}(\n\cap\k^\perp)} n_\beta \beta.
$$
In addition, $\omega = \mu -2\rho^\perp_\n$ by hypothesis.  Hence
$$
w(\mu + \rho ) - (\mu + \rho ) = \sum\limits_{\beta\in\supp_\t (\n\cap\k^\perp)} n_\beta \beta.
$$
\enddemo

On the other hand, since $\mu + \rho$ is $\b_\k$-dominant, there exists a multiset \newline $m$.: $\supp_\t (\n\cap\k)\to\N$ such that $(\mu + \rho ) - w(\mu + \rho ) = \sum\limits_{\alpha\in\supp_\t (\n\cap \k)} m_\alpha \alpha$.  Therefore
$$
\sum\limits_{\alpha\in\supp_\t (\n\cap\k )} m_\alpha \alpha + \sum\limits_{\beta\in\supp_\t (\n\cap \k^\perp )} n_\beta \beta = 0
$$
and
$$\sum\limits_{\alpha\in\supp_\t (\n\cap\k )} m_\alpha \langle\Re\lambda ,\alpha\rangle + \sum\limits_{\beta\in\supp_\t (\n\cap \k^\perp)} n_\beta \langle\Re\lambda, \beta\rangle = 0.
$$
Hence $m_\alpha = n_\beta = 0$ for all $\alpha, \beta$, and $w(\mu + \rho ) = \mu + \rho$.  As $\mu + \rho$ is a regular weight of  $\k$, $w= \text{id}$ and $i=0$. \qed

\proclaim{Lemma 2} $\dim\Hom_\k (V, F^s(\p,E)) = \dim E$.
\endproclaim

\demo{Proof}  Lemma 1 enables us to rewrite (1) in the special case $\delta = \mu$ as
$$\align
&\dim \Hom_\k (V(\mu), F^{s}(\p,E)) \\
&=\sum\limits_{0\leq j\leq s} (-1)^j (\sum\limits^\infty_{m=0} \dim \Hom_\t (H^j(\n\cap\k, V(\mu)),
S^m(\n\cap\k^\perp)\otimes E\otimes \Lambda^{\dim (\n\cap\k^{\perp})} (\n\cap\k^\perp))),
\endalign
$$
and, by Kostant's theorem, $\supp_\t H^\cdot (\n\cap\k, V(\mu)) =  \{\tilde\sigma (\mu +\rho )-\rho\mid\tilde \sigma\in W_\k\}$ and $\mu$ appears with multiplicity $1$ in $\{\tilde\sigma (\mu +\rho)-\rho\mid\tilde\sigma\in W_\k\}$.  On the other hand,
$$\align
&\supp_\t (S^\cdot (\n\cap\k^\perp)\otimes E\otimes \Lambda^{\dim (\n\cap\k^{\perp})} (\n\cap\k^\perp)) \\
&\qquad  = \{\mu + \sum\limits_{\beta\in\supp_\t (\n\cap\k^{\perp})} n_\beta ~\beta \mid n_\beta\in\N\}.
\endalign
$$
Since $\mu + \rho$ is $\b_\k$-dominant,
$$
\{\tilde \sigma (\mu +\rho)-\rho\mid\tilde\sigma\in W_\k \}\subset \{ \mu -\sum\limits_{\alpha\in\supp_\t (\n\cap\k)} m_\alpha \alpha\mid m_\alpha\in\N\}.
$$
This, together with the inequality $\langle \Re \lambda, \gamma\rangle > 0 ~ \forall\gamma\in\supp_{\frak t} \frak n$ (see the proof of Lemma 1), allows us to conclude that
$$
\{\tilde\sigma (\mu +\rho )-\rho\mid\tilde\sigma\in  W_\k\}\cap\{\mu + \sum\limits_{\beta\in\supp_\t (\n\cap\k^\perp)} n_\beta \beta\} = \{\mu \}.
$$
Consequently,
$$
\Hom_\t (H^j (\n\cap\k, V(\mu )), S^m (\n\cap\k^\perp)\otimes E\otimes \Lambda^{\dim (\n\cap\k^\perp)} (\n\cap\k^\perp ))\ne 0
$$
only for $m=0$.  This shows that
$$\align
&\dim\Hom_\k (V(\mu ), F^s (\p,E))  \\
&= \dim \Hom_\t (H^0 (\n\cap\k, V(\mu)), E\otimes \Lambda^{\dim(\n\cap\k^\perp)} (\n\cap\k^\perp)) = \dim E. \qed
\endalign
$$
\enddemo

\proclaim{Lemma 3}   If $V(\delta)$ is a $\k$-type of $F^s(\p,E)$ and $\delta\neq\mu$, then $\Vert\Re\delta + 2\rho\Vert > \Vert\Re\mu + 2\rho\Vert$.
\endproclaim

\demo{Proof}  By Proposition 1, and there exists a multiset $n$.: $\supp_\t (\n\cap \k^\perp) \to \N$ such that
$$
\delta + \rho = \mu + \rho + \sum\limits_{\beta\in\supp_\t (\n\cap\k^\perp)} n_\beta \beta.
$$
Hence
$$
\delta + 2\rho = \mu + 2\rho + \sum\limits_{\beta\in\supp_\t (\n\cap\k^\perp )}n_\beta \beta .
$$
Since $\p\subset\p_{\mu +2\rho}, \langle\Re\mu + 2\rho, \beta\rangle\geq 0$ for all $\beta\in\supp_\t (\n\cap\k^\perp )$.  In addition, $\delta\neq\mu$ implies $\Vert\sum\limits_{\beta\in\supp_\t (\n\cap\k^\perp )} n_\beta \beta\Vert^2 > 0$.
Therefore
$$\align
&\Vert\Re\delta + 2\rho\Vert^2 = \Vert\Re\mu +2\rho\Vert^2 + \Vert\sum\limits_{\beta\in\supp_\t (\n\cap\k^\perp )} n_\beta \beta\Vert^2 \\
&+ 2 \sum\limits_{\beta\in\supp_\t (\n\cap\k^\perp )} n_\beta \langle\Re\mu + 2\rho, \beta\rangle > \Vert\Re\mu + 2\rho\Vert^2. \qed
\endalign
$$
\enddemo

\head 3.  Discussion \endhead

An ultimate goal of the program of  study laid out in \cite{PZ1} is
the classification of simple generalized Harish-Chandra modules.
Within this framework, Theorem 1 above establishes the non-emptiness
of the class of simple $(\g,\k)$-modules of finite type with a fixed
minimal $\k$-type $V$, where $V$ is an arbitrary $\k$-type.  If $V =
V(\mu )$ is a generic $\k$-type (the definition, see \cite{PZ2},
involves certain inequalities on $\mu$), all modules in this class
are classified in \cite{PZ2} and in particular are subquotients of
$F^s(\p, E)$ generated by the unique minimal $\k$-type $V$ of
$F^s(\p,E)$ constructed exactly as in the present note as
subquotients of $F^s(\p,E)$ generated by $V$. For a non-generic $V$,
Theorem 2 yields an interesting class of simple generalized
Harish-Chandra modules which deserves further study.  It is known
that in general, these modules do not  exhaust all simple
generalized Harish-Chandra modules, as when the pair $(\g,\k )$ is
symmetric, or when $\k$ is a Cartan subalgebra of $\g$, the
classifications of simple $(\g,\k )$-modules in these two cases
yield modules which do not arise through our construction. For
instance, in the latter case no cuspidal modules, i.e., modules on
which all root vectors act freely, are fundamental series modules.
On the other hand, there are symmetric pairs $(\g,\k )$ for which
our construction yields all simple Harish-Chandra modules.  This
applies in particular to pairs of the form $(\s\oplus\s,\s )$, where
$\s$ is a simple Lie algebra and the inclusion
$\s\hookrightarrow\s\oplus\s$ is the diagonal map.  It is an
interesting question whether for some general (non-symmetric) pairs
$(\g,\k )$ the construction of this paper exhausts all simple
$(\g,\k )$-modules of finite type.

\vskip .20in

\noindent Ivan Penkov

\noindent International University Bremen

\noindent Campus Ring 1, D-28759

\noindent Bremen, Germany

\noindent email:  i.penkov\@iu-bremen.de
\vskip .20in

\noindent Gregg Zuckerman

\noindent Department of Mathemaics

\noindent Yale University

\noindent 10 Hillhouse Avenue

\noindent P.O. Box 208283

\noindent New Haven, CT 06520-8283, USA

\noindent email:  gregg\@math.yale.edu

\end